%% file: 2017_paperVRSSMUM.tex
\documentclass[[a4paper,1p]{elsarticle}

\usepackage{lineno,hyperref}
\modulolinenumbers[5]

\def\statement{\begin{minipage}[t]{.75\textwidth}
       NOTICE: this is the author's version of a work that was accepted for publication in Computers and Fluids. Changes resulting from the
       publishing process, such as peer review, editing, corrections, structural formatting, and other quality control mechanisms may not be reflected in this
       document. Changes may have been made to this work since it was submitted for publication. A definitive version was subsequently published in Computers and Fluids,
       \url{https://doi.org/10.1016/j.compfluid.2018.04.006} \\\\
       \textcopyright \ 2018. This manuscript version is made available under the CC-BY-NC-ND 4.0 license \url{http://creativecommons.org/licenses/by-nc-nd/4.0/}
       \end{minipage}}

\makeatletter
\def\ps@pprintTitle{%
     \let\@oddhead\@empty
     \let\@evenhead\@empty
     \def\@oddfoot{\footnotesize\itshape
       \statement\hfill\today}%
     \let\@evenfoot\@oddfoot}
\makeatother

\journal{Computers and Fluids}

\usepackage[font=small]{subcaption}
\usepackage{amsmath}
\usepackage{amsfonts}
\usepackage{booktabs}
\usepackage[section]{placeins}

\begin{document}

\begin{frontmatter}

\title{The Virtual Ring Shear-Slip Mesh Update Method}


\author{Fabian Key\corref{mycorrespondingauthor}}
\ead{key@cats.rwth-aachen.de}
\author{Lutz Pauli}
\author{Stefanie Elgeti}

\cortext[mycorrespondingauthor]{Corresponding author}

\address{Chair for Computational Analysis of Technical Systems, \\ CCES, RWTH Aachen University, \\ 52056 Aachen, Germany}




\begin{abstract}
    A novel method -- the Virtual Ring Shear-Slip Mesh Update Method (VR-SSMUM) -- for the efficient and accurate modeling of moving boundary or interface problems
    in the context of the numerical analysis of fluid flow is presented. We focus on cases with periodic straight-line translation including object entry and exit.
    The periodic character of the motion is reflected in the method via a mapping of the physical domain onto a closed virtual ring. Therefore, we use an extended
    mesh, where unneeded portions are deactivated to control the computational overhead. We provide a validation case as well as
    examples for the applicability of the method to 2D and 3D models of packaging machines.
\end{abstract}


\end{frontmatter}


\captionsetup[subfigure]{font=footnotesize,labelfont=footnotesize}

\section{Introduction}
    \input{src/introduction/introduction.tex}

\section{Problem Statement}
\label{sec:problemStatement}
        \input{src/problemStatement/problemStatementVRSSMUM.tex}

\section{Governing Equations and Solution Method}
\label{sec:govEqnsSolutionMethods}
    \input{src/implementation/spaceTimeFEM.tex}

\section{Mesh Update Method}
\label{sec:MUM}
    \input{src/implementation/implementationVRSSMUM.tex}

\section{Numerical Examples}
\label{sec:numericalExamples}
    \input{src/numericalExamples/numericalExamplesVRSSMUM.tex}

\section{Conclusion}
\label{sec:conclusion}
    \input{src/conclusion/conclusion.tex}

\section*{Acknowledgements}
The authors gratefully acknowledge the research funding, which was provided by SIG Combibloc Systems GmbH.
The computations were conducted on computing clusters supplied by the J{\"u}lich Aachen Research Alliance (JARA) and the RWTH IT Center.

\section*{References}

\bibliography{2017_paperVRSSMUM}

\end{document}

%% file: src/introduction/introduction.tex
\begin{figure}[h]
    \includegraphics[width=\textwidth]{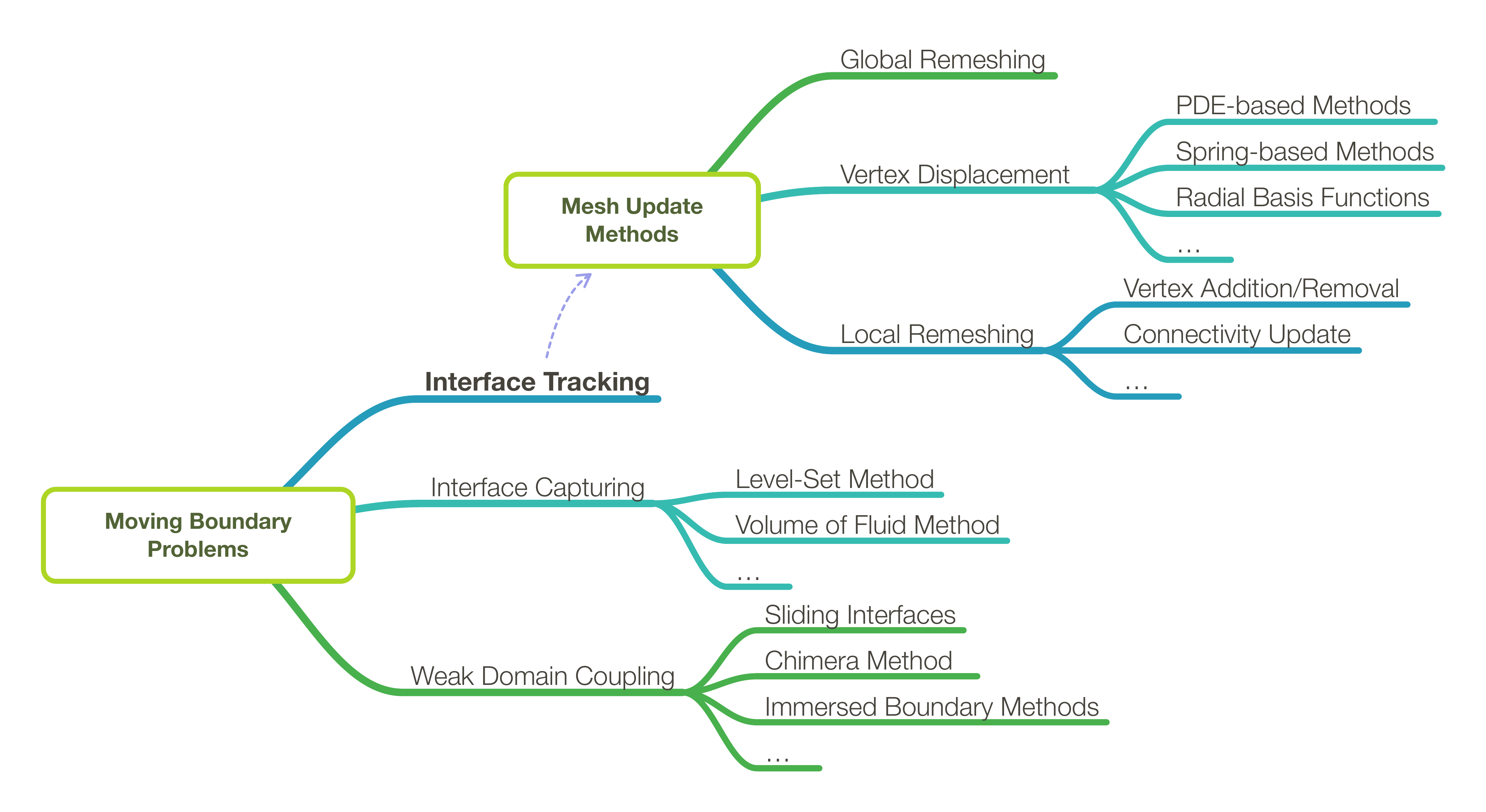}
    \caption{Overview of numerical methods for moving boundary problems.}
    \label{fig:mindMapMovingBoundProbs}
\end{figure}

As part of the analysis of fluid flow problems, we often encounter situations in which we have to deal with deforming domains. These can either be due to movement of outer boundaries or internal interfaces.
Problems including moving boundaries or interfaces impose special requirements on the numerical methods to model the involved motion with respect to both the computational mesh and the solution field.
Figure \ref{fig:mindMapMovingBoundProbs} shows an overview of different approaches that have been developed for the analysis of moving boundary or interface problems.
\\\\
Commonly, a first distinction is made between interface capturing and interface tracking methods \cite{elgeti2016}.
In the former approach, one employs an implicit description of the dynamic boundary or interface on some reference mesh.
More precisely, the distribution of a scalar function is used to identify the bounds of the considered domains.
This strategy implies a certain degree of flexibility, since (1) it allows complex interfaces and even topological changes, and (2) it operates on fixed spatial domains.
However, it can potentially involve drawbacks with respect to accuracy of mass conservation, resolution of both discontinuities and high gradients across the interface, or the imposition of boundary conditions.
The most common examples for interface capturing methods are the level-set method \cite{osher1988,chang1996} and the volume-of-fluid method \cite{hirt1981}.
\\\\
The class of interface tracking methods is based on boundary-conforming meshes. This requires an update of the computational mesh to account for the moving boundary.
There exist several strategies to adjust the mesh according to the change in the position of the interface.
One elementary approach is to perform a global remeshing of the altered geometry. As a consequence, the solution field has to be projected from the old mesh to the newly created one.
This increases not only the computational effort, but also introduces an additional projection error. In general, the computational costs for a complete remeshing are quite high.
Therefore, this is not a feasible solution for most cases.
\\
In contrast, other mesh deformation methods follow the idea of retaining at least the connectivity information of the original mesh.
They aim at varying only the position of interior mesh nodes based on the prescribed boundary motion; an approach, which reduces the computational effort significantly.
Examples are PDE-based methods like the elastic mesh update method \cite{tezduyar1992}, spring-based methods \cite{batina1990} or the concept of radial basis functions \cite{deBoer2007}.
\\
Some further methods enhance these strategies and apply a subsequent mesh optimization based on local mesh topology changes \cite{alauzet2013,wang2015}. These could be, e.g., edge swapping or vertex smoothing operations.
As a consequence, also large displacements can be considered.
\\\\
For situations in which the previous methods are not applicable, local remeshing strategies have been developed.
For cases including large relative movement, the Shear-Slip Mesh Update Method (SSMUM) was introduced and applied in \cite{tezduyar1996} to translational and in \cite{behr1999_SSMUM,behr2001,behr2003} to rotational movement.
Here, only elements in a small portion of the mesh are deformed and remeshed by means of a connectivity update. This circumvents a projection of the solution field and reduces the effort of updating the mesh.
So far, the method was mainly applied to model the relative movement of rotating objects with respect to other fixed parts. Several results for simulations in both 2D and 3D as well as on unstructured grids are available,
but there is only one example which covers translational movement.
\\
In the context of finite volume schemes, boundary movement in two dimensional problems has been modeled by adding and removing vertices in the corresponding space-time mesh \cite{zwart1998}.
\\\\
As an alternative to the former methods, one can also follow weak domain coupling strategies to account for the moving boundary or interface.
Typically, these methods are based on composite grids that are made up of sub-grids, which are not connected in the classical sense.
The coupling of the domains is rather achieved by introducing additional conditions on the solution field.
The Chimera method \cite{steger1983, benek1986, steger1987} was introduced to simplify the mesh generation for complex geometries. It works on overlapping grids that, individually, are easy to generate.
The coupling is performed via an interpolation of the solution field in the overlapping regions, which again requires a projection.
The continuity of the solution field can also be weakly imposed over a sliding interface, e.g., between rotating components \cite{bazilevs2008, takizawa2015}.
A further approach is the immersed boundary method \cite{peskin1972}, where simulations are always performed on a cartesian background grid. The actual boundary of the geometry is immersed in these grids and boundary conditions are applied
via an additional forcing function, which enters the underlying state equations.
\\\\
Summarizing, a great variety of methods exists for the solution of moving boundary problems. The choice of a specific method commonly depends on the special characteristics of the present problem.
Thus, we give a characterization of problems which we aim to solve with the proposed method in the following.
We will focus on cases that include large relative translational movement and assume periodic motion including topological domain changes in the sense of object entry and exit. Our aim is to avoid projection of the solution and the associated inaccuracies, as well
as global remeshing or mesh modifications to keep the computational effort small.
In this paper, we will present a novel method that allows to efficiently and accurately handle moving boundary problems. We developed an interface tracking approach that is embedded in a Deforming-Spatial Domain/Stabilized Space-Time (DSD/SST)
finite element framework. It is an extension of the SSMUM and its main idea is to map the translational movement in the physical domain
to a continuous circular movement in an abstract space, i.e., along a virtual ring. Therefore, we will refer to it as the Virtual Ring Shear-Slip Mesh Update Method (VR-SSMUM).
\\\\
The paper starts with a characterization of the underlying moving boundary problem and explains the fundamental ideas of the proposed mesh update method in Section \ref{sec:problemStatement}.
In Section \ref{sec:govEqnsSolutionMethods}, we will briefly outline the space-time finite element method for the incompressible Navier-Stokes equations
and place the focus on the implementation of the VR-SSMUM into this framework in Section \ref{sec:MUM}. Section \ref{sec:numericalExamples} contains numerical examples including a validation case and further test cases in 2D and 3D to demonstrate the applicability of the method.

%% file: src/problemStatement/problemStatementVRSSMUM.tex
\begin{figure}[h]
    \includegraphics[width=\textwidth]{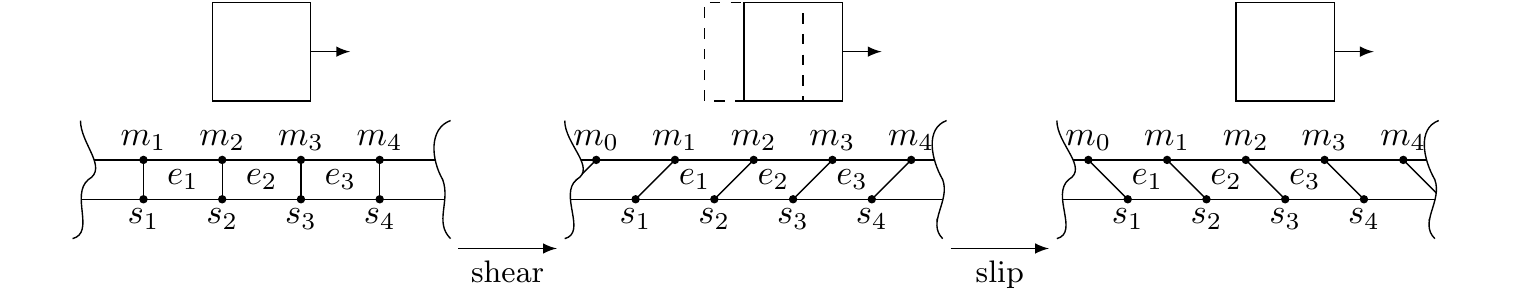}
    \caption
    {
        Elements $e_{\bullet}$ in the update layer deform due to movement of the square (\textit{shear} step). A connectivity update reverts the deformation (\textit{slip} step).
        Nodes of the moving and static mesh are denoted as $m_{\bullet}$ and $s_{\bullet}$, respectively.
    }
    \label{fig:movingBoxShearSlip}
\end{figure}
\begin{figure}
    \centering
    \begin{subfigure}[t]{0.49\textwidth}
        \includegraphics[width=\textwidth]{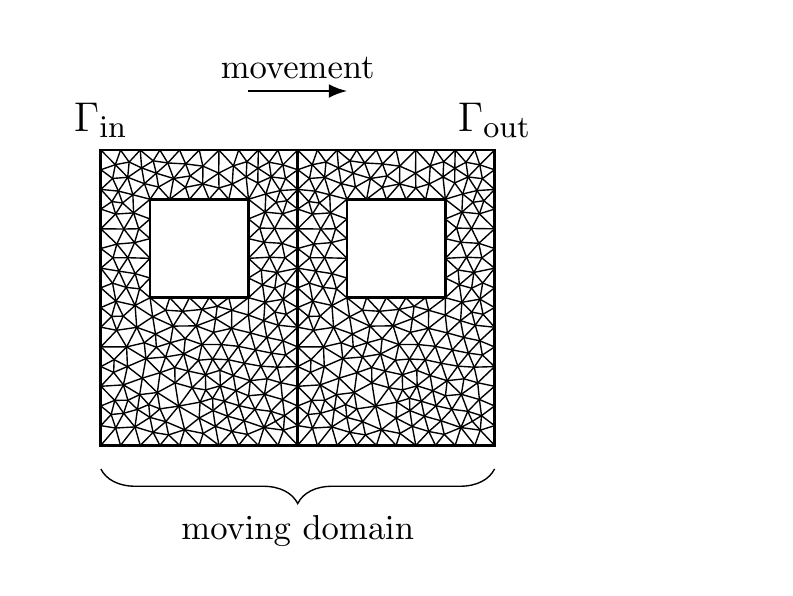}
        \caption{}
        \label{fig:movingCharacteristicBlocks}
    \end{subfigure}
    \centering
    \begin{subfigure}[t]{0.49\textwidth}
        \includegraphics[width=\textwidth]{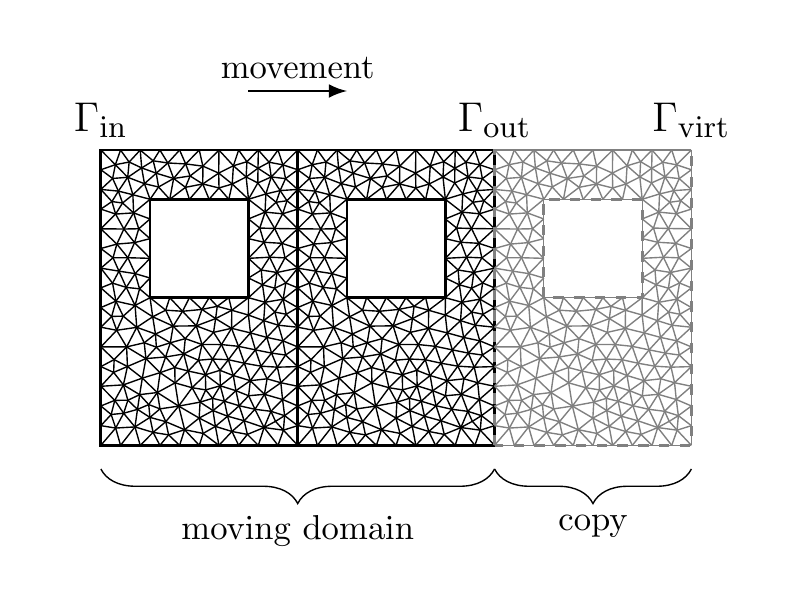}
        \caption{}
        \label{fig:movingCharacteristicBlocksPlusCopy}
    \end{subfigure}
    \caption
    {
        Moving domain for two squares which undergo translation.
        (a) Domain composed of characteristic mesh blocks.
        (b) Characteristic mesh blocks plus additional copy.
    }
    \label{fig:twoMovingSquares}
\end{figure}
\begin{figure}
    \centering
    \begin{subfigure}[t]{0.4\textwidth}
        \includegraphics[width=\textwidth]{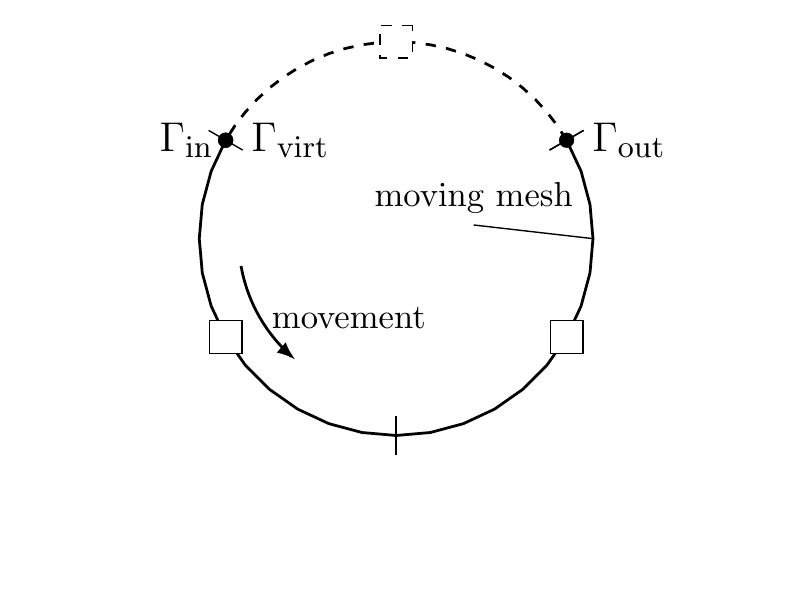}
        \caption{}
        \label{fig:movingVirtualRingClosed}
    \end{subfigure}
    \centering
    \begin{subfigure}[t]{0.4\textwidth}
        \includegraphics[width=\textwidth]{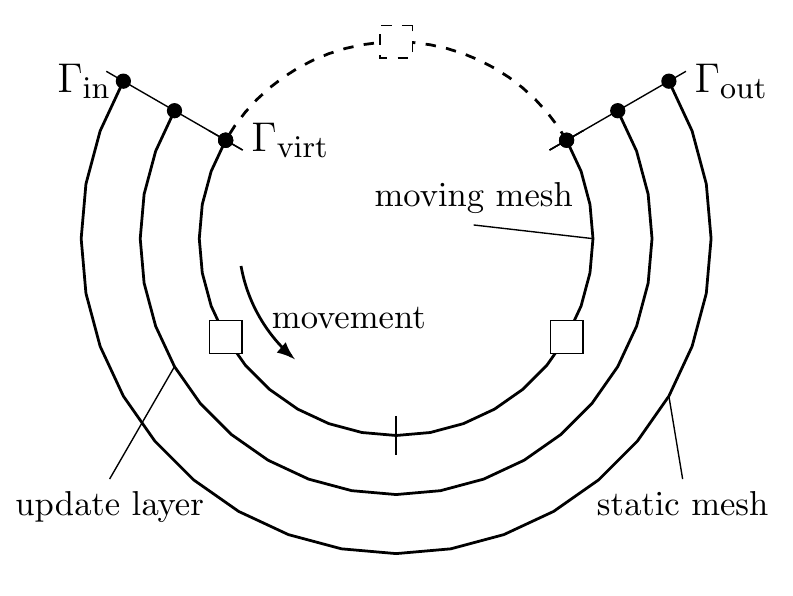}
        \caption{}
        \label{fig:movingVirtualRingClosedLinked}
    \end{subfigure}
    \caption
    {
        Virtual ring in abstract space for the example of two moving squares.
        (a) Moving mesh closed by additional mesh block (dashed).
        (b) Closed moving mesh combined with static mesh by update layer.
    }
\end{figure}

In the following, we will focus on the modeling of domain deformations that are on the one hand strictly translational, prescribed and periodic in nature and
on the other hand affect only a portion of the domain boundary, thus leading to large relative movement.
The Shear-Slip Mesh Update Method (SSMUM) presented in \cite{behr1999_SSMUM} has been developed for this kind of problems including large but regular boundary displacements.
Although the concept was proposed also for straight-line translation, only one example was shown in \cite{tezduyar1996}. Therein, two trains passing in a tunnel have been simulated.
The computations have been performed on a structured grid and the issue of modeling the entry or exit of the trains into or out of the tunnel has not been adressed.
\\\\
However, the applications we have in mind require the use of unstructured grids due to their complex geometry and also include object entry and exit by the virtue of the periodic movement.
As an example, this could be the simulation of processes which include the movement of a conveyor belt through a stationary machine casing as present in packaging machines or coating procedures.
Thus, we present the Virtual Ring Shear-Slip Mesh Update Method (VR-SSMUM), which is an extension of the SSMUM, and is applicable to a certain class of problems including translational movement. More specifically, the movement can be characterized in the way described above:
unidirectional, periodic and relative to a static reference boundary.
The method is able to work on unstructured grids and, additionally, allows us to model object entry and exit. The main benefits of the SSMUM, namely no need for either global remeshing or projection of the solution, are inherited.
This implies the favorable properties of the method with respect to efficiency and accuracy.
\begin{figure}
    \centering
    \begin{subfigure}[t]{0.4\textwidth}
        \includegraphics[width=\textwidth]{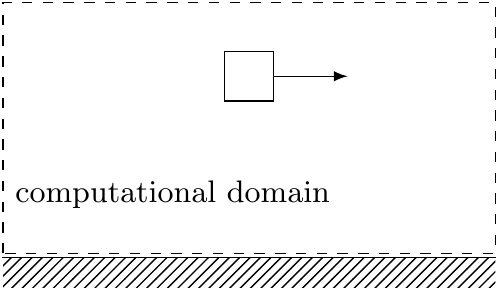}
        \caption{}
        \label{fig:movingBoxWall}
    \end{subfigure}
    \begin{subfigure}[t]{0.4\textwidth}
        \includegraphics[width=\textwidth]{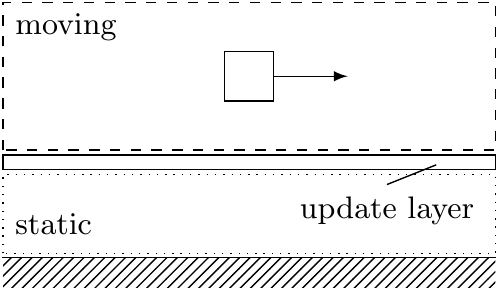}
        \caption{}
        \label{fig:movingBoxWallSplit}
    \end{subfigure}
    \caption
    {
        Two objects in relative motion: the square undergoes an unidirectional translation, whereas the wall does not move.
        (a) Computational domain.
        (b) Splitting into a moving and a static part connected by update layer.
    }
\end{figure}
\\\\
First, we will shortly review the concept of the SSMUM which served as basis for the presented VR-SSMUM in the context of translational movement.
The basic idea is to split the computational domain into a moving and a static part, respectively. Subsequently, a thin layer of elements -- the update layer -- is added, which connects the moving and static mesh portions.
For illustration, we will consider two objects in relative translational movement, e.g., a box which moves along a static wall (see Figure \ref{fig:movingBoxWall}).
\\\\
For this example, a splitting of the computational domain into a moving portion, a static portion and the update layer is shown in Figure \ref{fig:movingBoxWallSplit}.
The moving and static domains each hold a mesh, which is fixed with respect to the corresponding boundaries. This means that the moving mesh will perform a rigid body displacement as soon as the object and, thus, its boundary
starts to move.
Due to the movement of the interface between the moving mesh and the update layer, the elements in this layer undergo a shear deformation. This is referred to as \textit{shear} step in the SSMUM.
When the movement proceeds, the elements will become more and more distorted.
To counter this, the SSMUM applies a \textit{slip} step, which means that the update layer is remeshed by means of a connectivity update, i.e., the elements will then be comprised of a different, yet neighbouring set of nodes.
In Figure \ref{fig:movingBoxShearSlip}, this procedure is illustrated for the discussed example of the moving square.
\\
The connectivity update is based on update node information: each node of an affected element knows by which node it has to be replaced in the update procedure.
A projection of the nodal solution coefficients is not necessary, since the nodes retain their position.
Nevertheless, an update of the connectivity alters the basis functions in the affected elements even if the nodes stay in place.
Thus, the representation of the solution differs when the nodal values remain unchanged.
A modification of the original SSMUM circumvents this inaccuracy by changing the connectivity over a time slab instead of updating it between two time slabs \cite{schippke2012}.
Another approach is to incorporate the old connectivity during the integration of the jump term inside the update layer.
However, our experience is that the solutions with and without this additional treatment show no significant differences if the grid, i.e., the update layer, is fine enough.
\\\\
Next, we will outline the underlying idea and the additional features of the VR-SSMUM, which allows us to handle the class of moving boundary problems described above.
In order to address the periodicity, the moving domain is assumed to be built-up of characteristic blocks.
Furthermore, all blocks have identical discretization, i.e., the same mesh. In the most general case, our moving domain can consist of only one block.
Hence, the assumption of characteristic blocks does not entail any restrictions on the choice of the moving geometry. The setup for an example with two moving squares is shown in Figure \ref{fig:movingCharacteristicBlocks}.
\\\\
The modeling of the movement by means of the computational grid is achieved through the following idea: the translational, unidirectional movement in physical space is mapped to a movement along a virtual ring in a
more abstract space.
To set up this ring, the mesh for the moving domain is closed between its physical boundaries $\Gamma_{\text{in}}$ and $\Gamma_{\text{out}}$ by an additional copy of the characteristic mesh block (see Figure \ref{fig:movingCharacteristicBlocksPlusCopy}).
The boundary of this block will be denoted as $\Gamma_{\text{virt}}$. We let $\Gamma_{\text{in}}$ and $\Gamma_{\text{virt}}$ coincide in the abstract space and, thus, the moving mesh now represents a closed ring. In Figure \ref{fig:movingVirtualRingClosed},
the virtual ring in the abstract space is illustrated for the example discussed.
We can identify the original physical domain boundaries as points on this ring. Between these two points, the additional virtual block is depicted. The other segments correspond to the original blocks
of the moving domain.
As one can conclude from the figure, the geometric periodicity is automatically implied, since we have assumed that all blocks are identical. However, the geometric periodicity does not enforce periodicity with respect to
the solution field, since the mesh block, which enters the domain, is not the same one, which leaves the domain. For the solution process, the mesh portion in the virtual region between outlet and inlet is deactivated.
It is only used to model the mesh deformation.
\\\\
Once the movement can be modeled following this concept, we have to connect the moving and static mesh portions. This can be done in the same way as it has been proposed in the presentation of the SSMUM:
in the intermediate space between the different mesh portions, we add a thin layer of hand made elements, which aggregate the individual sub-meshes into one (see Figure \ref{fig:movingVirtualRingClosedLinked}).
\\\\
In contrast to the original version of the SSMUM, the update criterion is not (directly) based on the shape of the elements in the update layer. Instead, an update is triggered if any mesh point, which is part of the interface between
the moving mesh portion and the update layer, enters the virtual region and disappears in the physical space. Due to the identity of the mesh blocks, the corresponding copy of this point in the next block
concurrently enters the physical domain at the opposite boundary.
Thus, it is ensured that we are always able to perform the connectivity update in the update layer.

%% file: src/implementation/spaceTimeFEM.tex
In the present study, we deal with the flow field of a viscous, incompressible fluid. Let $t \in \left( 0,T \right)$ be an instant of time and $n_{\text{sd}}$ the number of space dimensions.
The time-dependent computational domain for time $t$ given by $\Omega_t \subset \mathbb{R}^{n_{\text{sd}}}$ is enclosed by its boundary $\Gamma_t$.
The velocity and pressure fields, $\mathbf{u}(\mathbf{x},t)$ and $p(\mathbf{x},t)$, evolve according to the incompressible Navier-Stokes equations.
The equations for the conservation of momentum and mass can be written as:
\begin{align}
	\rho
    \left(
        \frac{\partial \mathbf{u}}{\partial t} + \mathbf{u} \cdot \nabla \mathbf{u} - \mathbf{f}
    \right)
    - \mathbf{\nabla} \cdot \boldsymbol{\sigma}
    = \mathbf{0}
    \text{ on } \Omega_t \ \ \forall t \in (0,T),
    \label{eq:fluidNavierStokesMomentum}
    \\
	\mathbf{\nabla} \cdot \mathbf{u} = 0 \text{ on } \Omega_t \ \ \forall t \in (0,T),
    \label{eq:fluidNavierStokesMass}
\end{align}
where $\rho$ is the fluid density, which is assumed to be constant. An external force field enters the equations via the source term $\mathbf{f}(\mathbf{x},t)$.
Considering a Newtonian fluid, the following relation for the stress tensor $\boldsymbol{\sigma}$ is used to close the set of equations:
\begin{align}
    \boldsymbol{\sigma} \left( \mathbf{u},p\right) = - p \mathbf{I} + 2 \mu \boldsymbol{\varepsilon}(\mathbf{u}),
    \label{eq:stressTensor}
\end{align}
where
\begin{align}
\boldsymbol{\varepsilon}(\mathbf{u})
&=
\frac{1}{2}
\left(
\nabla \mathbf{u}
+ \left(
\nabla \mathbf{u}
\right)^T
\right).
\label{eq:rateOfStrainTensor}
\end{align}
The boundary conditions of Dirichlet or Neumann type read
\begin{align}
    \mathbf{u} = \mathbf{g} \text{ on } (\Gamma_t)_g, \\
    \mathbf{n} \cdot \boldsymbol{\sigma} = \mathbf{h} \text{ on } (\Gamma_t)_h,
\end{align}
with complementary portions $(\Gamma_t)_g$ and $(\Gamma_t)_h$ of $\Gamma_t$.
\\\\
The computations presented here are based on the DSD/SST finite element formulation \cite{tezduyar1992_DSD/SST_1}. Instead of constructing the weak form for the underlying equations only over the spatial domain, the space-time formulation involves the
corresponding space-time domain.
This framework is particularly suitable for moving boundary or interface problems, since the deformation of the spatial domain over time is naturally involved. In the discrete sense, the shape of the space-time elements
already implies the movement.
\begin{figure}
		\centering
    \includegraphics{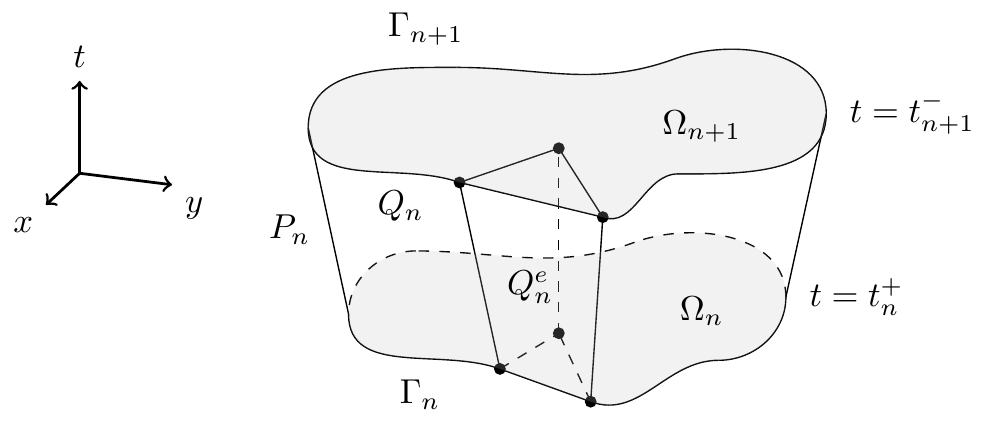}
    \caption{Sketch of a space-time slab $Q_n$ between lower and upper time level $t_n^+$ and $t_{n+1}^-$ with exemplary space-time element $Q_n^e$.}
    \label{fig:spaceTimeSlab}
\end{figure}
\\\\
Next, we will present the finite element function spaces for the DSD/SST method. Consider the time interval $(0,T)$ and a partition into subintervals $I_n = (t_n, t_{n+1})$.
The two time levels $t_n$ and $t_{n+1}$ line up in the ordered sequence of time levels $ 0 = t_0 < t_1 < ... < t_N = T$.
For the time level $n$, the spatial domain and its boundary are denoted as $\Omega_n = \Omega_{t_n}$ and $\Gamma_n = \Gamma_{t_n}$, respectively.
A so called space-time slab $Q_n$  is defined as the domain that emerges between $\Omega_n$, $\Omega_{n+1}$ and the surface $P_n$, which emerges from $\Gamma_t$ when $t$ passes $I_n$.
Figure \ref{fig:spaceTimeSlab} shows a sketch of such a space-time slab. Similar to the spatial boundary $\Gamma_t$, we can assign portions of $P_n$, $(P_n)_g$ and  $(P_n)_h$, belonging to the space-time boundary of Dirichlet or Neumann type.
\\\\
For each space-time slab, the finite element function spaces for first-order polynomials in space and time are given as
\begin{align}
(\mathbf{\mathcal{S}}^h_{\mathbf{u}})_n &= \lbrace \mathbf{u}^h \in \lbrack H^{1h}(Q_n) \rbrack^{n_{\text{sd}}} \ | \ \mathbf{u}^h = \mathbf{g} \text{ on } (P_n)_g \rbrace,
\\
(\mathbf{\mathcal{V}}^h_{\mathbf{u}})_n   &= \lbrace \mathbf{w}^h \in \lbrack H^{1h}(Q_n) \rbrack^{n_{\text{sd}}} \ | \ \mathbf{w}^h = \mathbf{0} \text{ on } (P_n)_g \rbrace,
\\
(\mathcal{S}^h_p)_n  &= (\mathcal{V}^h_p)_n  = \lbrace p^h \in H^{1h}(Q_n)  \rbrace.
\end{align}

In contrast to the original publications, we follow the slightly modified formulation presented in \cite{pauli2017}.
The stabilized space-time formulation for the incompressible Navier-Stokes equations (\ref{eq:fluidNavierStokesMomentum}) and (\ref{eq:fluidNavierStokesMass}) then reads as follows:
\\\\
\textit{Given $(\mathbf{u}^h)^-_n$, find $\mathbf{u}^h \in (\mathbf{\mathcal{S}}^h_{\mathbf{u}})_n $ and $p^h\in (\mathcal{S}^h_p)_n $ such that
        $\forall \mathbf{w}^h \in (\mathbf{\mathcal{V}}^h_{\mathbf{u}})_n $ and $\forall q^h \in (\mathcal{V}^h_p)_n   $:}
\begin{align}
    \nonumber
    \int_{Q_n} {\bf w}^h \cdot \rho \Big( \frac{\partial {\bf u}^h}{\partial t} + {{\bf u}^h \cdot \boldsymbol{{\nabla}} {\bf u}^h} - {\bf f} \Big)dQ
    + \int_{Q_n} \mbox{\boldmath$\nabla$}{\bf w}^h:\mbox{\boldmath$\sigma$} ({\bf u}^h, p^h)dQ \\ \nonumber
    + \int_{Q_n} q^h \boldsymbol{{\nabla}} \cdot {\bf u}^hdQ
    + \int_{\Omega_n} ({\bf w}^h)^+_n \cdot \rho \left( \left({\bf u}^h\right)^+_n - \left({\bf u}^h\right)^-_n \right) d\Omega\\ \nonumber
    + \sum_{e=1}^{({n_{\text{el}})}_n} \int_{Q^e_n} \tau_{\mbox{\tiny{MOM}}} \frac{1}{\rho}
    \Big[ \rho \Big( {{\bf u}^h \cdot \boldsymbol{{\nabla}} {\bf w}^h} \Big) + \boldsymbol{{\nabla}} q^h \Big] \\ \nonumber
    \cdot \Big[ \rho \Big(\frac{\partial {\bf u}^h}{\partial t} + {{\bf u}^h \cdot \boldsymbol{{\nabla}}  {\bf u}^h} - {\bf f} \Big) - \boldsymbol{{\nabla}} \cdot \mbox{\boldmath$\sigma$}({\bf u}^h, p^h) \Big] dQ \\
    + \sum_{e=1}^{({n_{\text{el}})}_n} \int_{Q^e_n} \tau_{\mbox{\tiny{CONT}}} \boldsymbol{{\nabla}}\cdot {\bf w}^h \rho \boldsymbol{{\nabla}}\cdot {\bf u}^h dQ
    = \int_{(P_n)_h} {\bf w}^h \cdot {\bf h}^h dP.
    \label{eq:stabWeakFormNavierStokes}
\end{align}
Here, we make use of the following notation:
\begin{align}
(\mathbf{u}^h)^\pm_n = \lim\limits_{\zeta \rightarrow 0} \mathbf{u}(t_n\pm \zeta)
\\
\int_{Q_n} \ldots dQ = \int_{I_n} \int_{\Omega_t} \ldots d\Omega dt
\\
\int_{P_n} \ldots dP = \int_{I_n} \int_{\Gamma_t} \ldots d\Gamma dt
\end{align}
The fourth term in Equation \ref{eq:stabWeakFormNavierStokes} is the so-called \textit{jump term}.
The purpose of this integral is to induce a weak continuity in time for the velocity field over $\Omega_n$ -- which is $\left({\bf u}^h\right)^+_n$ -- with respect to the previous space-time slab.
More precisely, we minimize the difference between $\left({\bf u}^h\right)^+_n $ and the velocity field over $\Omega_n$ from the underlying slab $\left({\bf u}^h\right)^-_n$.
For the stabilization parameters $\tau_{\mbox{\tiny{MOM}}}$ and $\tau_{\mbox{\tiny{CONT}}}$, we use the expressions given in \cite{pauli2017}.

%% file: src/implementation/implementationVRSSMUM.tex
\begin{figure}[h]
    \centering
    \includegraphics[width=\textwidth]{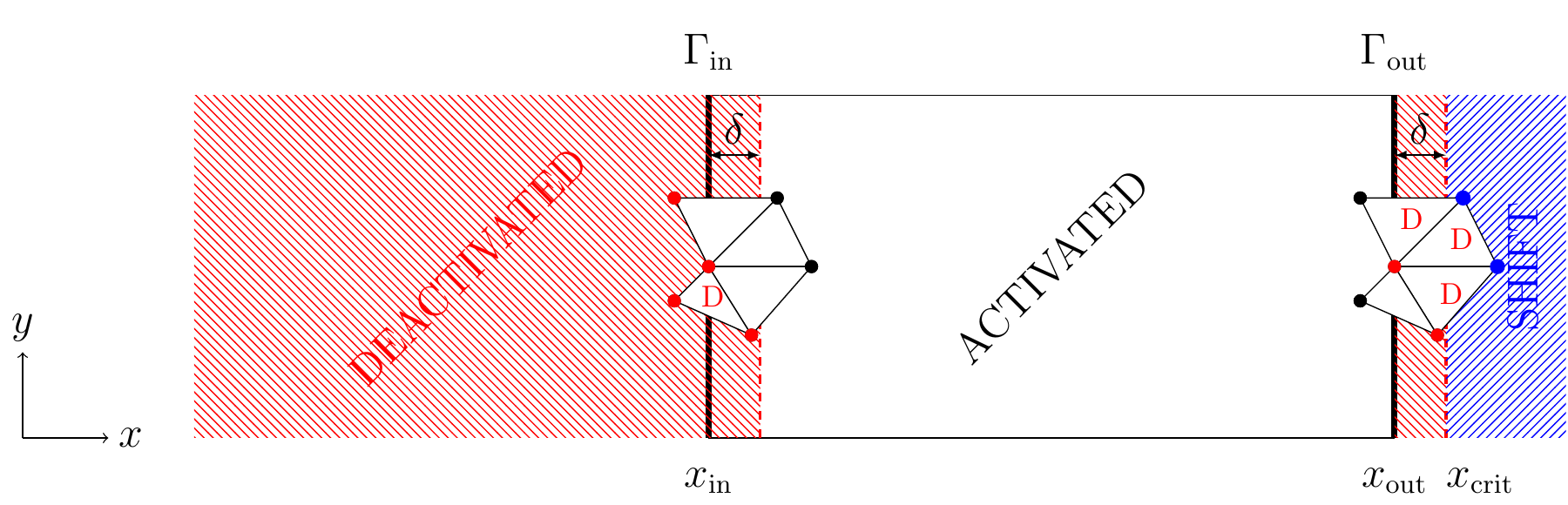}
    \caption
    {
        Moving domain with critical coordinate $x_{\text{crit}}$.
        Additionally, the zones for activating, deactivating and shifting nodes are depicted.
        As an example, multiple spatial elements are sketched. Their nodes are activated (black), deactivated (red) or shifted (blue).
        Elements with an ``D'' are deactivated, since they consist either of deactivated nodes only or hold a shifted one.
    }
    \label{fig:zones}
\end{figure}
\begin{figure}[h]
    \centering
    \includegraphics[width=\textwidth]{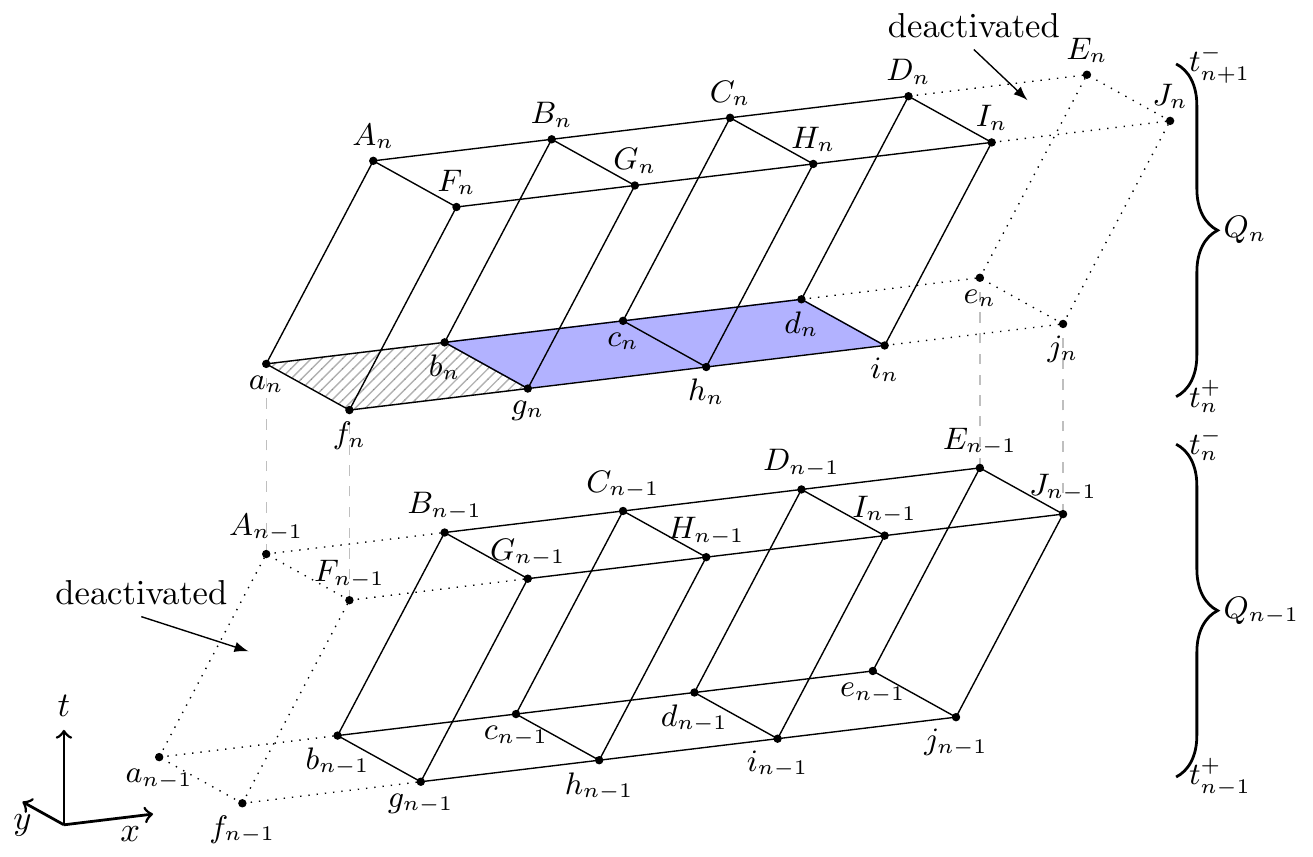}
    \caption
    {
        Example for changing element activity pattern and non-matching domains for two successive space-time slabs.
        Space-time elements with dotted edges are deactivated.
        The jump term integration is only performed over the solid-colored spatial elements.
        The shaded one is skipped, since we have no old solution for the nodes $a_n$ and $f_n$.
        Note that each space-time element has 8 nodes.
    }
    \label{fig:jumpTermNewActiveElementNotFixed}
\end{figure}
\begin{figure}[h]
    \centering
    \begin{subfigure}[t]{0.4\textwidth}
        \includegraphics[width=\textwidth]{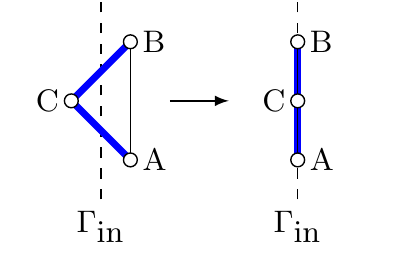}
        \caption{}
        \label{fig:collapsingElement}
    \end{subfigure}
    \begin{subfigure}[t]{0.4\textwidth}
        \includegraphics[width=\textwidth]{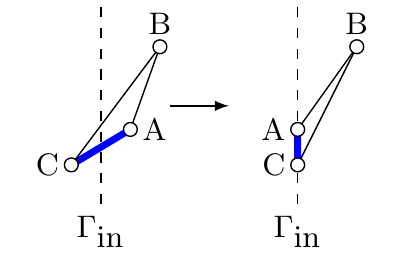}
        \caption{}
        \label{fig:twistedElement}
    \end{subfigure}
    \caption
    {
        Two examples where enforced boundary conformity leads to invalid spatial elements of an unstructured grid:
        (a) collapsing element, (b) twisted element.
        Boundary faces are drawn as thick blue lines.
    }
    \label{fig:invalidElements}
\end{figure}

In this section, we illustrate the implementation of the VR-SSMUM into the framework of the space-time FEM as presented in Section \ref{sec:govEqnsSolutionMethods}.
The VR-SSMUM mainly tackles two problems: first, it allows us to deal with large relative movement between moving and static objects (``SSMUM part'').
Additionally, the concept of the virtual ring makes it possible to model the object entry and exit (``VR part'').
Note that this concept -- taken individually -- may also be transferred to a semi-discrete approach, i.e., an Arbitrary Lagrangian Eulerian (ALE) formulation, with appropriate modifications.
The space-time formulation simply provided a framework in which the concept of the virtual ring could be conveniently combined with the SSMUM to accomplish the requirements of the applications discussed.
We will first outline the mesh generation process that is necessary to produce the special purpose computational grids which are suitable for the use with the VR-SSMUM.
Afterwards, the implications for the implementation of the space-time FEM are stated.
\\\\
Following the division of the computational domain into moving and static portions, an individual mesh is created for each part. The part of the domain which is allocated for the update layer remains empty.
The mesh generation for the static part does not require special treatment except for the discretization of the interface adjacent to the update layer. Details are presented below.
For the moving mesh, we have assumed in Section \ref{sec:problemStatement} that the corresponding domain is composed of several characteristic blocks. Thus, we start with a mesh for the first block and copy this sample to fill up the moving domain.
Additionally, one further copy is attached at the end of the original domain. In the next step, this copy will be used to close the ring.
As discussed in Section \ref{sec:problemStatement}, the closure happens through concurrence of the lateral boundaries $\Gamma_{\text{in}}$ and $\Gamma_{\text{virt}}$.
To achieve this also in the discrete sense, the indices of nodes on $\Gamma_{\text{virt}}$ are replaced by the corresponding indices of nodes on $\Gamma_{\text{in}}$. In the abstract sense, the mesh is closed now in the direction of movement.
In physical space, we have elements that span from one end to the other.
\\\\
The interface between the update layer and the moving or static mesh will be denoted as $\Gamma_{\text{M}}$ and $\Gamma_{\text{S}}$, respectively. We assume that we first discretize opposite interfaces in a structured manner and, if necessary,
switch to an unstructured one afterwards by diagonalization.
As a restriction, the discretization of $\Gamma_{\text{M}}$ has to provide at least as many elements as $\Gamma_{\text{S}}$ along the direction of movement.
This ensures that we are able to create the elements in the update layer which connect the moving and static mesh.
Furthermore, the structure of the interface discretization allows us to assign for each node on $\Gamma_{\text{M}}$ an unique update node, i.e., the node which will replace the original node in case of a recreation of the connectivity in the update layer.
When we have set up the individual computational grids as described above, we can manually create the layer of elements which connects the static and moving parts.
\\\\
The modifications regarding the procedure of the FEM are presented hereinafter.
The mesh update process per time step can be broken down into the following steps:
\begin{enumerate}
    \item Perform the connectivity update if the corresponding criterion is fulfilled.
    \item Apply the prescribed movement and determine active nodes and elements.
    \item Update the lateral boundaries.
    \item Incorporate the activity information during assembly and solution of the linear system.
\end{enumerate}
At the beginning of each time step, we have to decide if a connectivity update in the update layer has to be performed for the mesh of the related space-time slab.
Strictly speaking, we have stated in Section \ref{sec:problemStatement} that a node is deactivated as soon as it crosses the boundary $\Gamma_{\text{out}}$.
However, this condition is slightly relaxed in the actual implementation.
Therefore, we define a critical coordinate $x_{\text{crit}}$ with respect to the direction of movement, which is located further downstream. If $x_{\text{out}}$ is the position of $\Gamma_{\text{out}}$,
we define
\begin{equation}
    x_{\text{crit}} = x_{\text{out}} + \delta,
\end{equation}
with some distance $\delta$ defining two margins at the lateral boundaries $\Gamma_{\text{in}}$ and $\Gamma_{\text{out}}$ (see Figure \ref{fig:zones}).
The following criterion is used to detect a connectivity update:
check if a node on the interface $\Gamma_{\text{M}}$ exceeds the critical coordinate $x_{\text{crit}}$.
If the criterion is fulfilled, we reconnect the corresponding elements in the update layer.
As already mentioned above, this procedure is based on the update node information.
\\\\
The next step is to update the nodal coordinates in the moving mesh according to the prescribed motion.
Since we follow the space-time approach outlined in Section \ref{sec:govEqnsSolutionMethods}, the prescribed movement can be directly applied to the space-time mesh of the current time step.
This is accomplished by altering the coordinates of the nodes on the upper time level. Therefore, the current velocity is evaluated and the resulting displacement is applied to the nodes on the upper time level.
Although the nodes can continuously move along the virtual ring, a shift oft their physical coordinates is necessary at some point. In other words, nodes which have already left the physical domain need to jump to its front such that they can enter
the domain again when the movement proceeds.
We apply this shift with respect to the physical coordinates to those nodes which have moved too far. Here, we again make use of the critical coordinate from the update criterion.
Every node which has left the right margin, i.e., it has crossed $x_{\text{crit}}$, has reached the \textit{shift region} and is shifted to the deactivation zone in front of the domain.
\\\\
Based on the updated nodal coordinates, we determine on the upper time level the nodes which remain in our computational domain and, thus, stay activated. Therefore, we use the (de-)activation zones depicted in Figure \ref{fig:zones}.
We always assign the same activity information to the associated nodes on the lower time level.
Using the nodal activity pattern, we set the activated elements. The following rule is applied: every element with at least one activated node is also activated; otherwise, it is deactivated.
Furthermore, all elements containing a node that has been shifted are deactivated.
Due to the setup with multiple identical mesh blocks, an element at the boundary $\Gamma_{\text{in}}$ has a related copy at $\Gamma_{\text{out}}$.
The aforementioned rule ensures that always only one instance of these elements is activated at a time. Furthermore, the distance $\delta$ should be less or equal the minimal extent of any element in the direction of movement.
This prevents that both elements are deactivated at the same time, since at least one of the element nodes is always located in the activation zone.
By now, boundary nodes have been deactivated to turn on or off the appropriate elements. To correct this and to take them into account in the computation, all nodes belonging to active elements are activated afterwards.
\\\\
Since the mesh on the virtual ring moves through the physical domain, we constantly have to update the element faces which form the boundaries $\Gamma_{\text{in}}$ and $\Gamma_{\text{out}}$.
The new lateral boundaries are identified as transition of activated to deactivated elements. Hence, also the nodes belonging to boundary faces change during the movement.
Since we allow to prescribe boundary conditions on $\Gamma_{\text{in}}$ and $\Gamma_{\text{out}}$, we have to update the active degrees of freedom before we start with the assembly of the linear system of equations.
Note that the VR-SSMUM itself does not restrict the type of boundary condition which is set on the lateral boundaries.
The choice of these boundary conditions rather depends on the specific problem which is investigated.
\\\\
During the assembly process, we also have to consider the element activity information. In general, all deactivated elements should not lead to a contribution to the system matrix and right hand side.
We achieve this goal by skipping the deactivated elements in the course of computing the element contributions.
As we will see in the following, the jump term from Equation (\ref{eq:stabWeakFormNavierStokes}) has to be treated carefully during the assembly process.
The integral of the jump term is defined over the spatial domain $\Omega_n$ at time $t_n$. Since we update the nodal coordinates and the element activity information for each space-time slab, new elements can enter the computational domain
compared to the previous one. As a result, the spatial domain, which is covered by the elements of the moving mesh, may be different for two successive space-time slabs.
According to the rule for activating elements, a newly activated element contains a number of nodes, which have been activated already in an earlier time step, and at least one node which just became active.
For the latter type of nodes, no information about the old solution is available.
Therefore, active elements, which have been deactivated in the previous time slab, are not considered in the integration of the jump term.
An example for two consecutive space-time slabs with different element activity patterns and, thus, non-matching spatial domains is shown in Figure \ref{fig:jumpTermNewActiveElementNotFixed}.
\\\\
A solution that would allow to incorporate newly activated elements in the jump term integration would need to provide information about the old solution for the previously deactivated nodes.
One idea would be the projection of the solution field from the inlet boundary surface of the previous space-time slab to the old solution field at the inlet of the current slab, which now consists of a different set of elements.
\\\\
As soon as the assembly process is finished, we can solve the resulting linear system. Here, we only include those degrees of freedom which are neither constrained by a boundary condition nor related to deactivated nodes.
Thus, the computational overhead arising from the additionally introduced mesh elements and nodes is limited.
\\\\
When the moving mesh undergoes a rigid body displacement, its boundaries do not conform with the original lateral boundaries of the computational domain anymore.
For structured grids, the boundary conformity of the moving portion can be achieved by relocating the nodes of the boundary faces onto the original domain boundary.
However, this approach cannot be easily adopted for unstructured grids, since it could lead to collapsing or twisted spatial elements. An example for each case is sketched in Figure \ref{fig:invalidElements}.
For the case shown in Figure \ref{fig:collapsingElement}, we can detect the collapsing element and ignore it for the current time step.
A solution for the case in Figure \ref{fig:twistedElement} would require a reordering of the element nodes.
Although this approach is not yet included in the current implementation, we observed that this situation only occurs very rarely.
Furthermore, if we can ensure that we do not encounter this situation for a specific simulation, no special treatment is necessary.
\\\\
Our mesh update method was integrated in a fully parallelized implementation of the DSD/SST method based on domain decomposition in use with MPI.
The integration did not require any special treatment with respect to parallelization.

%% file: src/numericalExamples/numericalExamplesVRSSMUM.tex
\subsection{Validation Case: 2D Couette Flow}
\begin{figure}
    \centering
    \begin{subfigure}{0.49\textwidth}
        \includegraphics[height=0.8\textwidth]{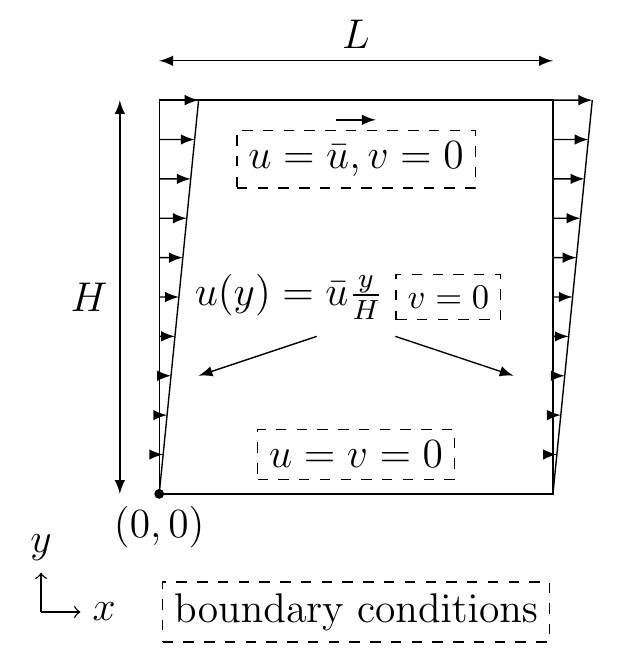}
        \caption{}
        \label{fig:couetteSetup}
    \end{subfigure}
    \begin{subfigure}{0.49\textwidth}
        \includegraphics[height=0.8\textwidth]{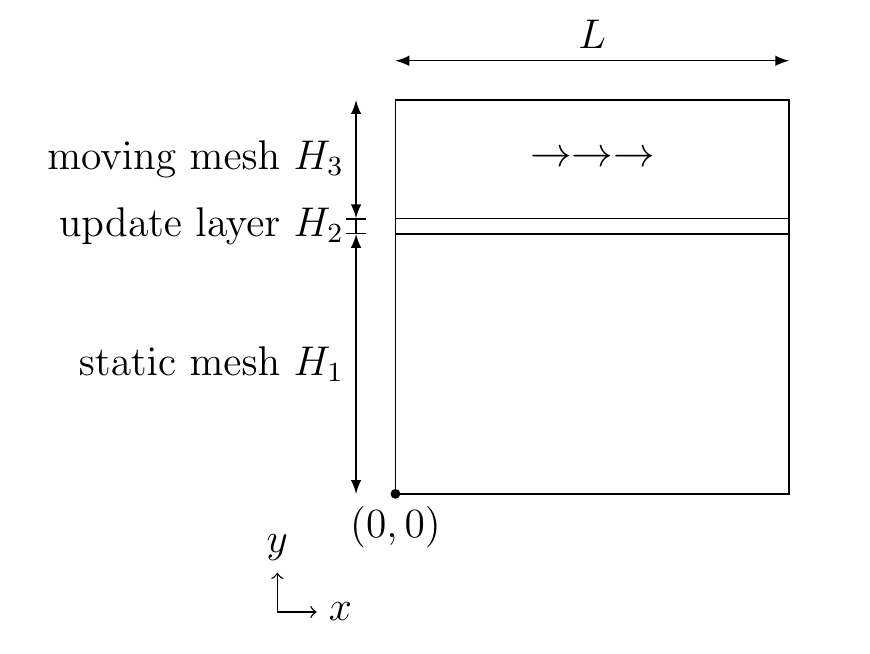}
        \caption{}
        \label{fig:couetteMesh}
    \end{subfigure}
    \caption
    {
        2D Couette flow:
        (a) setup with boundary conditions,
        (b) splitting of the computational domain into static and moving portions connected by the update layer.
    }
    \label{fig:couette}
\end{figure}

\begin{figure}
    \includegraphics[width=\textwidth]{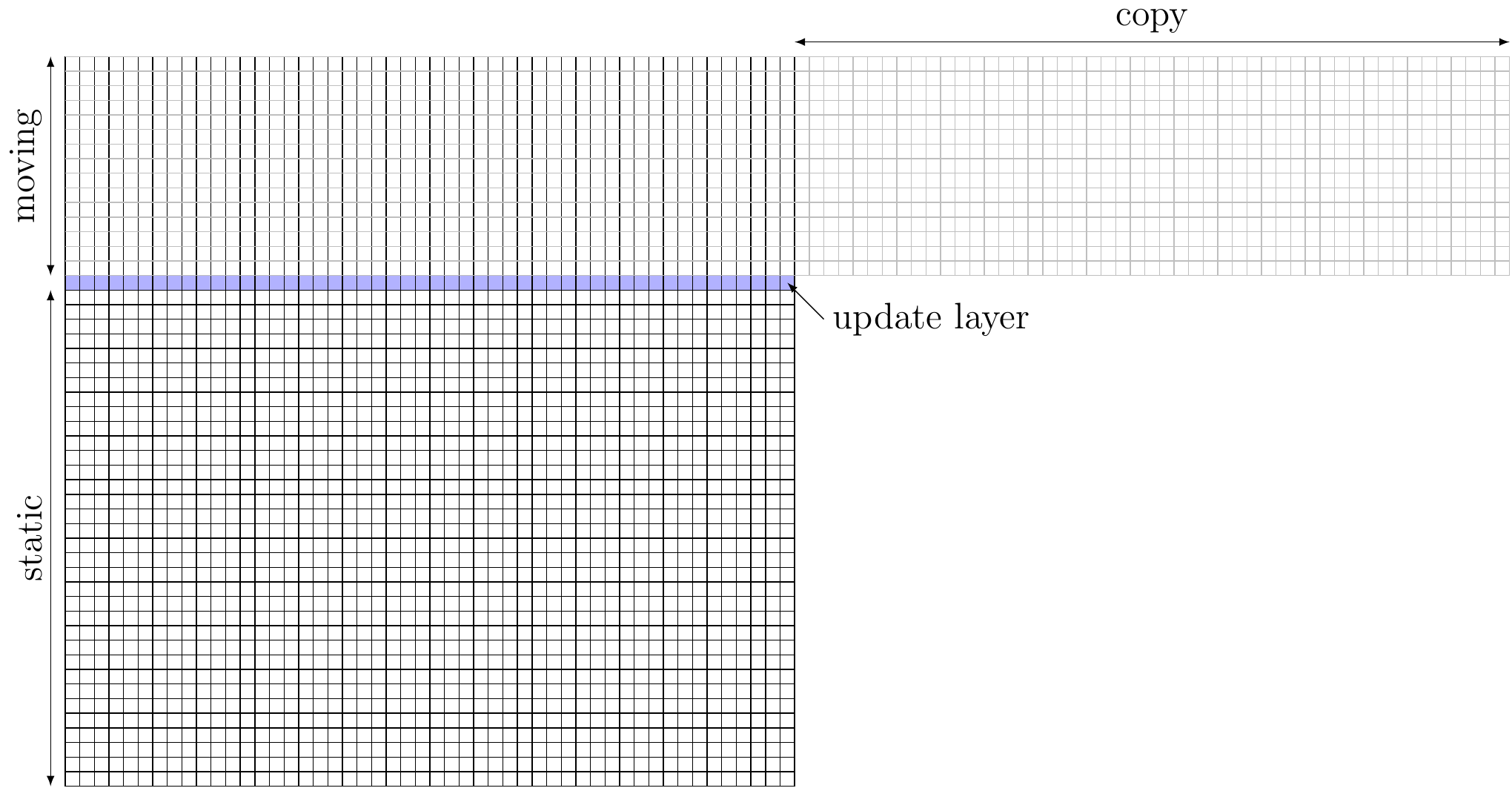}
    \caption{2D Couette flow: mesh prepared for VR-SSMUM.}
    \label{fig:couetteMeshVRSSMUM}
\end{figure}

\begin{figure}
    \centering
    \begin{subfigure}{0.49\textwidth}
        \includegraphics[width=\textwidth]{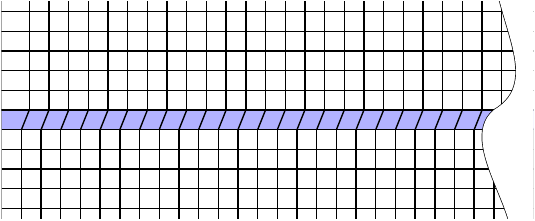}
        \caption{}
        \label{fig:couetteMeshUDL_5}
    \end{subfigure}
    \begin{subfigure}{0.49\textwidth}
        \includegraphics[width=\textwidth]{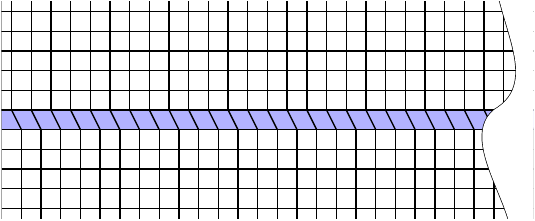}
        \caption{}
        \label{fig:couetteMeshUDL_6}
    \end{subfigure}
    \caption
    {
        2D Couette flow: part of the mesh with update layer (blue) for time steps 5 (left) and 6 (right).
        In time step 6, the element connectivity in this layer has been recreated to counteract
        the deformation resulting from the movement of the upper plate.
    }
    \label{fig:couetteMeshUDL}
\end{figure}

For the validation of the newly implemented VR-SSMUM, we make use of the classical Couette flow test case. The analytic solution, which is known in this situation, will allow us to confirm the correctness
of the numerical results. The setup of the test case is depicted in Figure \ref{fig:couetteSetup}. It naturally fits into the framework of moving and static objects treated by the VR-SSMUM. On the one hand, we have the fixed bottom plate.
On the other hand, we have the moving top plate. Together, we have two objects in relative, unidirectional and translational motion. The moving mesh portion will be connected to the upper plate, whereas the static
mesh portion contains the lower one. In between, the placement of the update layer is arbitrary. However, we place it off-center to avoid potential symmetry effects (see Figure \ref{fig:couetteMesh}).
Furthermore, we use a structured grid to illustrate also the conformity of the lateral boundaries. Figure \ref{fig:couetteMeshVRSSMUM} shows the mesh prepared for the VR-SSMUM.
It includes the additional copy for the moving domain, which consists of one original block in this case. The entire mesh contains 3,250 space-time elements with 6,770 nodes, where each moving block consists of 750 elements.
\\\\
For the setting of the Couette flow test case shown in Figure \ref{fig:couetteSetup}, the exact velocity distribution $\mathbf{u}=\left( u, v \right)^T$ is given at steady state as
\begin{align}
    u(y) &= \bar{u} \frac{y}{H}, \\
    v &\equiv 0,
\end{align}
where we have the velocity of the upper plate $\bar{u}$ and the distance between the two plates $H$. The dimensions and material parameters used for the simulation are given in Table \ref{tbl:2DCouetteFlowSettings}.
Information of the prepared mesh for the VR-SSMUM is stated in Table \ref{tbl:2DCouetteFlowMesh}. Since the dynamic mesh movement only allows unsteady simulations, the idea is to prescribe the analytic solution as initial conditions and to validate that the method does not distort the correct solution.
The relative error
\begin{align}
    \delta_{\text{rel}}^u = \frac{\left| \tilde{u}-u \right|}{\bar{u}}
\end{align}
with the numerical solution $\tilde{u}$ is evaluated for the first 8 time steps in Table \ref{tbl:relErrorU}.
\begin{table}
    \centering
    \begin{tabular}{cccccc}
        \toprule
        $H$ [m]& $L$ [m] & $\bar{u} \left[ \frac{\text{m}}{\text{s}} \right]$ & $\rho \left[\frac{\text{kg}}{\text{m}^3} \right]$ & $\mu \left[\frac{\text{kg}}{\text{m$\cdot$s}} \right]$\\
        1.0 & 1.0 & 0.02 & 100 & 2.5  \\
        \bottomrule
    \end{tabular}
    \caption
    {
        2D Couette flow: dimensions and parameter values.
    }
    \label{tbl:2DCouetteFlowSettings}
\end{table}
\begin{table}
    \centering
    \begin{tabular}{ccccc}
        \toprule
        $H_1$ [m]& $H_2$ [m]& $H_3$ [m] & $\Delta x$ [m]& $\Delta y$ [m]\\
        0.68 & 0.02 & 0.3 & 0.02 & 0.02 \\
        \bottomrule
    \end{tabular}
    \caption
    {
        2D Couette flow: mesh information.
    }
    \label{tbl:2DCouetteFlowMesh}
\end{table}
\begin{table}
    \centering
    \begin{tabular}{ccccc}
        \toprule
        Step & 1 - 5 & \textbf{6} & 7 & 8 \\
        $\delta_{\text{rel}}^u \cdot 10^{14}$ &
        0.0173 & 
        \textbf{5.98}   & 
        0.260  & 
        0.0520 \\ 
        \bottomrule
    \end{tabular}
    \caption
    {
        2D Couette flow:
        relative error for $x$-velocity $\delta_{\text{rel}}^u$.
        In time step 6, a connectivity update was performed and, more important,
        previously deactivated nodes became active.
    }
    \label{tbl:relErrorU}

\end{table}
For time steps 1 to 5, we cannot see any negative influence of our mesh update method. In time step 6, the relative error increases although it still accomplishes a good level of accuracy.
The reason for this effect is the fact that previously deactivated nodes become active at that moment. As already discussed in Section \ref{sec:MUM}, we have no information about the solution at these nodes from the previous time step.
This leads to a perturbation of the solution in the first non-linear iteration. However, we observe that the induced error vanishes over subsequent iterations.
To illustrate this situation, Figure \ref{fig:couetteMeshUDL} shows the mesh for time steps 5 and 6. For time step 5, the movement of the upper plate and the corresponding mesh portion
led to a deformation of the elements in the update layer. However, the original connectivity is still retained. In time step 6, new elements or rather nodes enter the domain at the left boundary accompanied by
a connectivity recreation of elements in the update layer.

\subsection{Test Case: 2D Packaging Machine with Object Entry and Exit}
\begin{figure}
    \centering
    \includegraphics[height=0.5\textwidth]{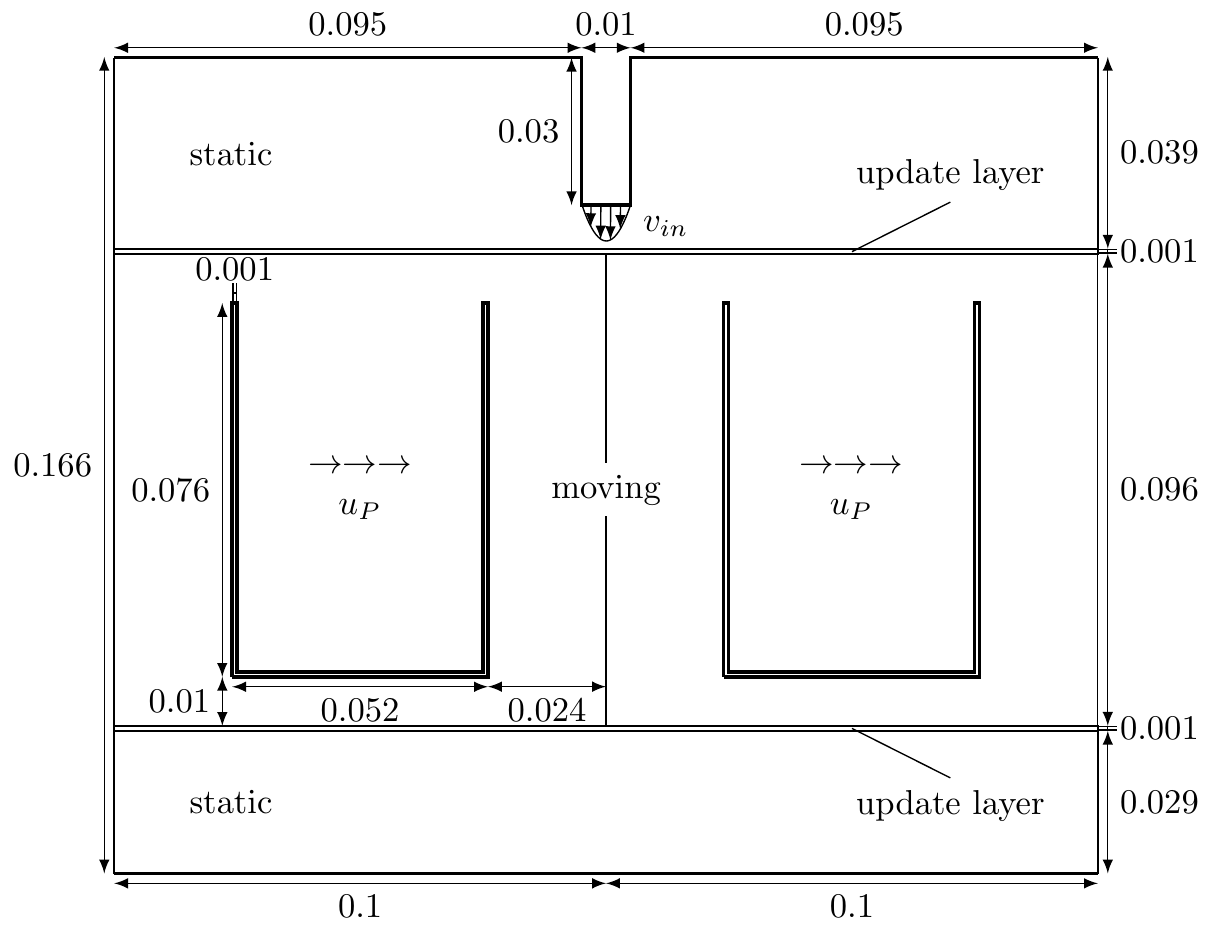}
    \caption
    {
        2D simplified packaging machine:
        setup and boundary conditions.
        No-slip walls are drawn as thick lines.
        All dimensions are given in meters.
    }
    \label{fig:2DsimplifiedPackagingMachineSetup}
\end{figure}

\begin{figure}
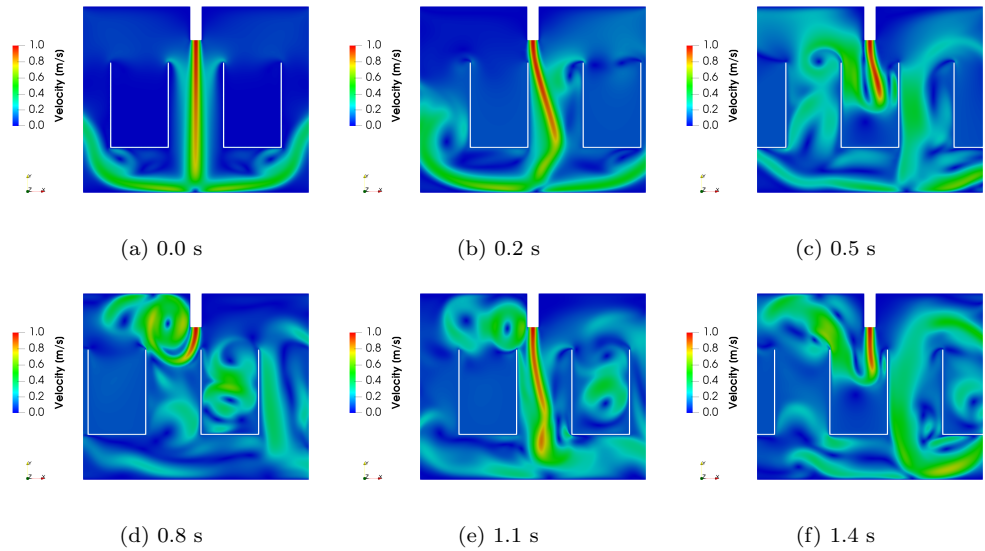

    \begin{subfigure}{0.32\textwidth}
        \includegraphics[width=\textwidth]{fig/numericalExamples/2DsimplifiedPackagingMachine/velocity_0.pdf}
        \caption{$0.0$ s}
        \label{}
    \end{subfigure}
    \begin{subfigure}{0.32\textwidth}
        \includegraphics[width=\textwidth]{fig/numericalExamples/2DsimplifiedPackagingMachine/velocity_100.pdf}
        \caption{$0.2$ s}
        \label{}
    \end{subfigure}
    \begin{subfigure}{0.32\textwidth}
        \includegraphics[width=\textwidth]{fig/numericalExamples/2DsimplifiedPackagingMachine/velocity_250.pdf}
        \caption{$0.5$ s}
        \label{}
    \end{subfigure}
    \\
    \begin{subfigure}{0.32\textwidth}
        \includegraphics[width=\textwidth]{fig/numericalExamples/2DsimplifiedPackagingMachine/velocity_400.pdf}
        \caption{$0.8$ s}
        \label{}
    \end{subfigure}
    \begin{subfigure}{0.32\textwidth}
        \includegraphics[width=\textwidth]{fig/numericalExamples/2DsimplifiedPackagingMachine/velocity_550.pdf}
        \caption{$1.1$ s}
        \label{}
    \end{subfigure}
    \begin{subfigure}{0.32\textwidth}
        \includegraphics[width=\textwidth]{fig/numericalExamples/2DsimplifiedPackagingMachine/velocity_700.pdf}
        \caption{$1.4$ s}
        \label{}
    \end{subfigure}
    \caption
    {
        2D simplified packaging machine: snapshots of the velocity field for several times.
        Packages pass the nozzle and enter or leave the domain.
    }
    \label{fig:2DsimplifiedPackagingMachine_velocity}
\end{figure}

In the following, we will show the application of the VR-SSMUM for a generic packaging machine.
We consider a simplified setup which consists of multiple packages passing through a machine casing.
At the top of the casing, a nozzle is located, which provides the inflow of fluid with velocity
\begin{align}
   v_{\text{in}}(x) = \bar{v} \cdot 40 \cdot 10^3(x-0.095)(x-0.105).
\end{align}
Figure \ref{fig:2DsimplifiedPackagingMachineSetup} shows a sketch of this setting.
Again, the computational domain is split into static and moving parts. In this case, we have two static mesh portions, i.e., the lower and the upper part,
which contain the machine casing and the nozzle. Consequently, also two update layers are present, which link the static mesh portions to the moving one, in
which the packages are embedded. For the mesh, the overall number of space-time elements is 154,970 with 158,210 nodes in total. The time step size is chosen
as $\Delta t = 2 \cdot 10^{-3}$ s and we run the simulation for 700 time steps.
\\\\
The movement of the packages is prescribed via a constant velocity $u_{\text{P}}$, which is set along the wall of the packages
to achieve a no-slip condition. The static walls of the machine casing are equipped with a no-slip condition as well.
Furthermore, the movement requires that we are able to let packages leave and enter the computational domain.
Although this type of boundary condition is typically suitable for boundaries far away from significant flow patterns,
we set zero Neumann conditions on the lateral boundaries, since this academic case is intended to rather illustrate the functionality of the method than to investigate the detailed flow field.
Note that we only consider the flow of a single-phase fluid. The parameter values for density
and viscosity are chosen to match hot air.
\\\\
We use the results for the static scenario with stationary packages ($u_{\text{P}}=0$), but injected fluid ($v_{\text{in}} \neq 0$), as initial condition.
Snapshots of the flow field at several points of time are given in Figure \ref{fig:2DsimplifiedPackagingMachine_velocity}.
Right after the beginning of the movement, we can observe how the packages carry along surrounding fluid and, that way, influence the inflow jet of the nozzle.
When a package passes the nozzle, the jet enters the inner space and induces several vortex structures. Afterwards, the jet leaves the package and the
procedure begins again for the subsequent package.
Note that we encounter situations in which the moving domain contains two complete packages as well as situations with one whole package and two packages,
which are only partially inside the domain.
\begin{table}
    \centering
    \begin{tabular}{ccccc}
        \toprule
        $u_{\text{P}} \left[ \frac{\text{m}}{\text{s}} \right]$& $\bar{v} \left[ \frac{\text{m}}{\text{s}} \right]$ & $\rho \left[\frac{\text{kg}}{\text{m}^3} \right]$ & $\mu \left[\frac{\text{kg}}{\text{m$\cdot$s}} \right]$\\
        0.1 & 1.0 & 0.6924 & 271 $\cdot 10^{-7}$ \\
        \bottomrule
    \end{tabular}
    \caption
    {
        2D simplified packaging machine: parameter values.
    }
    \label{tbl:2DsimplifiedPackagingMachineSettings}
\end{table}
\subsection{Use Case: 3D Packaging Machine}
\begin{figure}
    \begin{subfigure}{0.49\textwidth}
        \includegraphics[width=\textwidth]{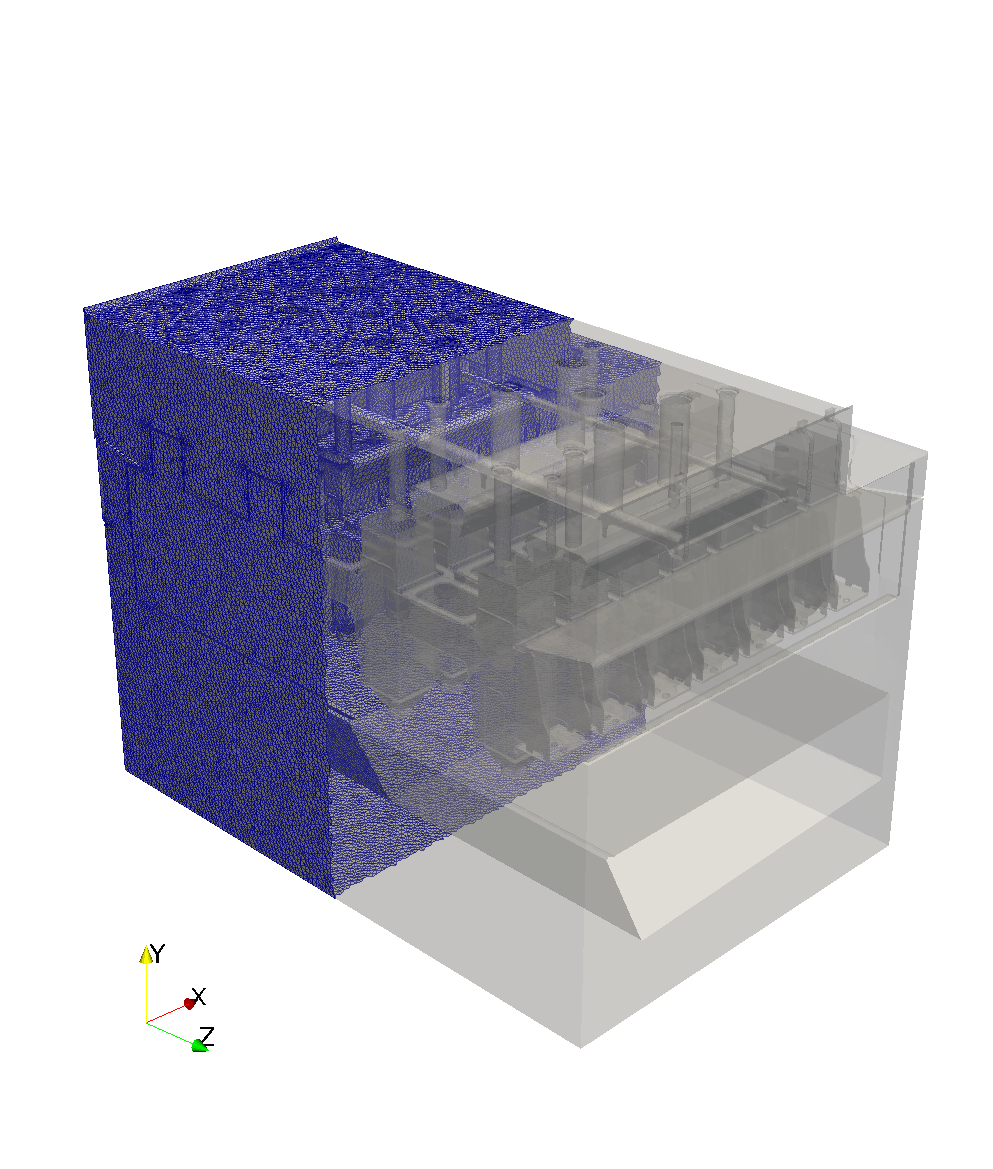}
        \caption{}
        \label{fig:3DpackagingMachineMesh}
    \end{subfigure}
    \begin{subfigure}{0.49\textwidth}
        \includegraphics[width=\textwidth]{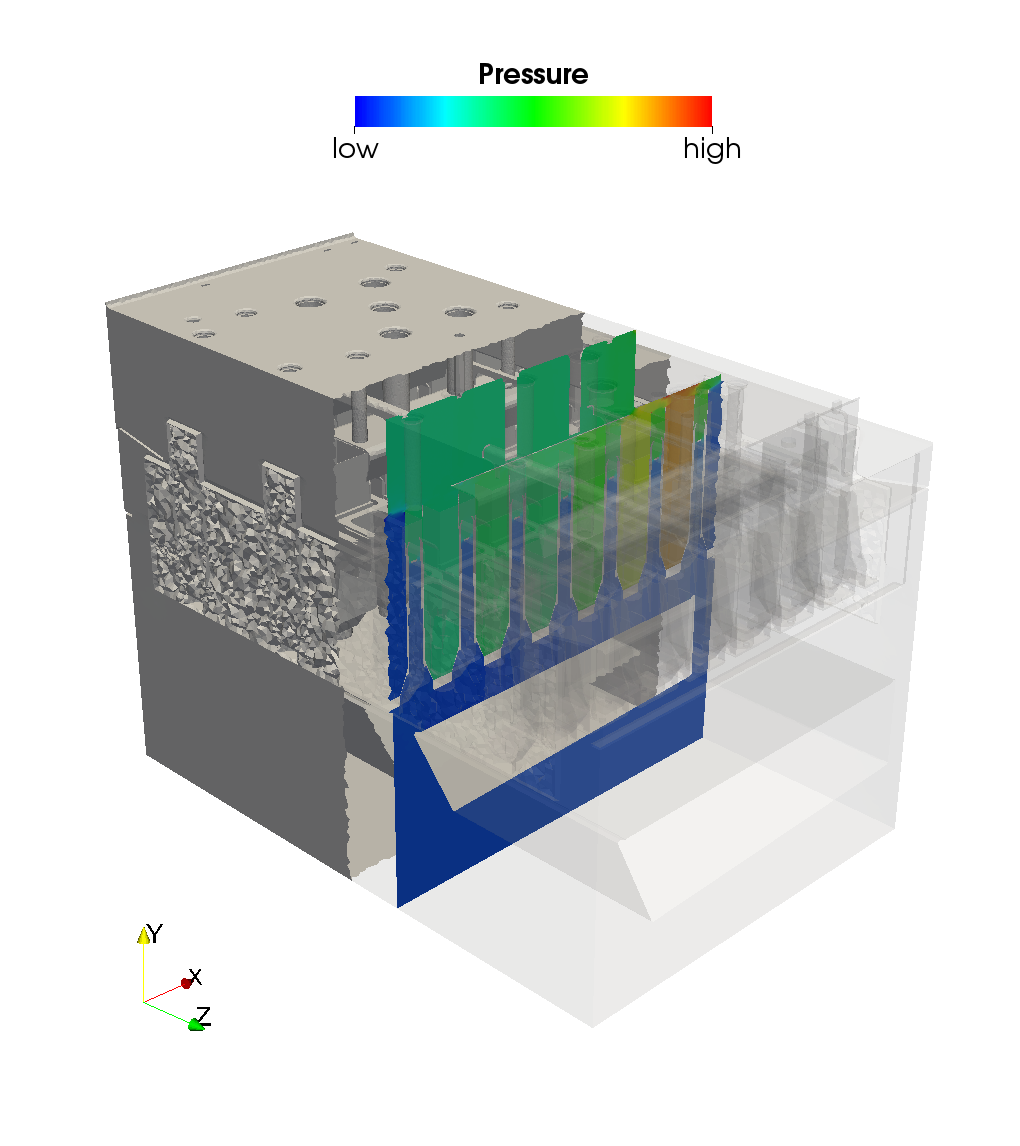}
        \caption{}
        \label{fig:3DpackagingMachinePressure}
    \end{subfigure}
    \caption
    {
        3D packaging machine:
        (a) mesh from real geometry data prepared for the VR-SSMUM,
        (b) slice with pressure distribution.
    }
    \label{fig:3DpackagingMachine}
\end{figure}

In this section, we will demonstrate that the presented mesh update method also works for industrial applications.
Therefore, we show results for a real three-dimensional packaging machine.
The basic setup of the machine is similar to the generic two-dimensional case which we discussed before.
Again, we have a static machine casing which contains several nozzles and further components at the top.
At the bottom, an additional suction box is mounted which removes redundant fluid.
Through the static machine casing, a conveyor belt transports the so-called pocket chain which holds the packages.
Based on the realistic geometry data, we prepared a mesh suitable for the VR-SSMUM. More precisely, we created a single
mesh for the static mesh portion, i.e., the outer machine casing with attached components, and a mesh for the moving block
containing the pocket chain and the packages. As described in Section \ref{sec:MUM}, we initially started with structured meshes for the interfaces,
which we subdivided afterwards to obtain complete unstructured meshes on the one hand but to make the update node information available on the other hand.
The resulting mesh is depicted in Figure \ref{fig:3DpackagingMachineMesh}. It is made up of 5,497,172 space-time elements and 2,198,944 nodes.
The pocket chain moves along the $x$-direction according to a prescribed time-dependent functional, which
describes acceleration and deceleration during one stroke of the machine.
In Figure \ref{fig:3DpackagingMachinePressure}, the pressure field is visualized on a slice through the machine after a certain movement of the
pocket chain.

%% file: src/conclusion/conclusion.tex
We presented the Virtual Ring Shear-Slip Mesh Update Method (VR-SSMUM) for the efficient and accurate modeling of
moving boundary or interface problems including straight-line translation, potentially on unstructured grids and in combination with object entry and exit.
As an extension of the Shear-Slip Mesh Update Method (SSMUM), the key idea is to map the physical movement onto a virtual ring in
an abstract space by extending the computational mesh properly.
The method was integrated into a framework based on the Deformable-Spatial-Domain/Stabilized Space-Time (DSD/SST) finite element formulation.
The favorable efficiency and accuracy properties of the SSMUM are automatically inherited. The computational overhead of the method is limited, since unneeded portions of the mesh are deactivated
during the solution process, i.e., the element assembly process and the solution of the linear system.
We validated the method by means of a Couette flow test case with analytic solution for the velocity field.
In addition, we showed the application to a relevant two-dimensional test case including object entry and exit based on a simplified packaging machine.
Finally, we performed a simulation for an industrial three-dimensional use case of a realistic packaging machine to prove the usefulness of the method in this field.